\documentclass{article}
\usepackage{amsmath}%
\usepackage{amsfonts}%
\usepackage{amssymb}%
\usepackage{color}%
\usepackage{graphicx}
\usepackage[normalem]{ulem}
\newtheorem{theorem}{Theorem}

\newtheorem{corollary}{Corollary}

\newtheorem{definition}{Definition}

\newtheorem{lemma}{Lemma}

\newtheorem{proposition}{Proposition}
\newtheorem{remark}{Remark}

\numberwithin{equation}{section}

\newcommand{\Mbb}{\mathbb}

\newcommand{\rr}{\mathbb R}

\newcommand{\rthr}{\mathbb R^3}
\newcommand{\Div}{\mathrm{div}}

\newcommand{\tr}{\mathrm{t}}
\newcommand{\OT}{\Omega_{T}}
\newcommand{\Oh}{\Omega_{h}}

\newcommand{\MK}{\mathbb{K}}
\newcommand{\vv}{\vec{v}}
\newcommand{\vn}{\vec{n}}
\newcommand{\vb}{\vec{b}}
\newcommand{\vw}{\vec{w}}

\newcommand{\Ned}{\mathcal{N}_{\epsilon+\delta}}

\newcommand{\Nepdl}{\mathcal{N}_{\eps,\delta}}

\newcommand{\Nepzr}{\mathcal{N}_{\eps,0}}
\newcommand{\Nepizr}{\mathcal{N}_{\eps_i,0}}

\newcommand{\zetaed}{\zeta_{\epsilon,\delta}}

\newcommand{\zetadd}{\zeta_{\delta,\delta}}
\newcommand{\zetaedt}{\partial_t\zeta_{\epsilon,\delta}}
\newcommand{\zetaddt}{\partial_t\zeta_{\delta,\delta}}

\newcommand{\aed}{a_{\epsilon,\delta}}
\newcommand{\aeidi}{a_{\epsilon_i,\delta_i}}

\newcommand{\aedt}{\partial_t a_{\epsilon,\delta}}

\newcommand{\ued}{u_{\epsilon,\delta}}
\newcommand{\ueidi}{u_{\epsilon_i,\delta_i}}

\newcommand{\eps}{\epsilon}

\newcommand{\ccu}{\mathfrak{u}}
\newcommand{\ccuC}{\mathfrak{u}_C}

\newcommand{\iOh}{\int_{\Omega_{h}}}

\newcommand{\JG}{\left|\partial_s\Gamma\right|}
\newcommand{\dGt}{d\delta_{\Gamma(t)}}

\newcommand{\vved}{\vec{v}_{\epsilon,\delta}}

\newcommand{\esup}{\mathrm{ess}\,\mathrm{sup}}

\newcommand{\bbox}{\vrule height.6em width.6em 
depth0em} 

\def\Xint#1{\mathchoice
	{\XXint\displaystyle\textstyle{#1}}%
	{\XXint\textstyle\scriptstyle{#1}}%
	{\XXint\scriptstyle\scriptscriptstyle{#1}}%
	{\XXint\scriptscriptstyle\scriptscriptstyle{#1}}%
	\!\int}
\def\XXint#1#2#3{{\setbox0=\hbox{$#1{#2#3}{\int}$ }
		\vcenter{\hbox{$#2#3$ }}\kern-.6\wd0}}

\def\dashint{\Xint-}

\let\TeXchi\chi
\newbox\chibox
\setbox0 \hbox{\mathsurround0pt $\TeXchi$}
\setbox\chibox \hbox{\raise\dp0 \box 0 }
\def\chi{\copy\chibox}

\begin{document}
\title{A Concentrated Capacity Model for Diffusion-Advection: \\ Advection Localized to a Moving Curve}
\author{Colin Klaus\\
             Mathematical Biosciences Institute, The Ohio State University\\  
			1735 Neil Avenue, Coumbus OH 43212, USA\\
email: {\tt klaus.68@osu.edu}}

\date{}
\maketitle
\vskip.4truecm
\begin{abstract}
In this work I show how a diffusion-advection equation in three space-dimensions may have its advection term weakly limited to a velocity field localized to a moving curve. This is rigorously accomplished through the technique of concentrated capacity, and the form of the concentrated capacity limit along with small time existence of solutions is determined. This problem is motivated by mathematical biology and the study of proteins in solvent where the latter is modeled as a diffusing quantity and the protein is taken to be a 1d object which advects the solvent by contact and its own motion. This work introduces a novel PDE's framework for that interaction.
\end{abstract}
\section{Introduction}\label{S:1}
Consider a bounded domain $\Omega\subset\rthr$ and its parabolic cylinder $\OT = \left[0,T\right)\times\Omega$.  Initially, we consider a diffusion-advection equation with vanishing Neumann condition
\begin{equation}
u_{t} + \Div\left[-\MK(x;t)\nabla u + u \vec{v} \right] = 0 \mbox{ in } \OT
\label{eq:diff_adv}
\end{equation}
\begin{equation}
\left(-\MK(x;t)\nabla u + u\vec{v}\right)\cdot \vec{n} = 0 \mbox{ on } \partial\OT
\label{eq:neu}
\end{equation}
We now would like to make precise the solution to such a problem when (\ref{eq:diff_adv},\ref{eq:neu}) is interpreted weakly and $\vec{v} = \vec{v}(t,s)d\delta_{\Gamma(t)}$ is Dirac-mass localized to a moving curve. However, as stated, such a formulation lacks the necessary compactness to justify the existence of solutions in a natural function space. The technical difficulty is its producing uncontrolled energy terms of the form $\int_{\Gamma(t)} u\nabla u \vec{v}ds$ after taking $u$ as a test function. It is natural to replace the distributional velocity field $\vec{v}$ with an approximation $\vec{v}_{\epsilon}$ which recovers the distributional field weakly in the limit. In this paper, it is demonstrated how the technique of \textit{concentrating capacity} makes this intuition precise.  

Other instances of concentrated capacity in the literature may be found in \cite{Showalter89,Colli90,Magenes95,Magenes96} for example. More recent examples are \cite{Andreucci03,Keller17}. A common theme is that a subset of the domain, which sites essential physics but is intrinsically lower-dimensional up to a small parameter $\eps$, is diminished by passing $\eps\rightarrow 0$. Sufficient compactness is included so that, in the limit, a new partial differential equation is recovered with the true lower dimensional domain and which preserves the physics of the original problem. 
\subsection{Novelty and Significance}
The investigation presented here is novel for its applying concentrating capacity to a time-varying curve rather than a stationary hypersurface. In particular, the technique is shown to make precise a parabolic PDE whose equation coefficients are highly singular, time-varying measures. 

In addition to its mathematical novelty, this problem takes inspiration from computational biology and the study of proteins in solvent. Proteins are chains of amino acids linked together in sequence along a polypeptide bond and are the molecular machines that drive many dynamics inside biological cells. In molecular dynamics simulations \cite{McCammon06,Shaw12}, resolving the motion of the explicit solvent which surrounds the protein is a great computational expense \cite[Sec 8.3]{McCammon06}. However, the localized interaction of solvent at different sites along the protein can strongly influence the conformation the protein takes \cite{Bellissent16,Laage17}, which in turn affects its ability to bind its targets or otherwise properly function. The implications for this extend to the molecular basis for disease as well as drug-design \cite{Amaro16b}. This concentrated capacity model presents a novel PDE's framework for considering these issues, while future investigations may incorporate other essential features such as amino acid sidechains and molecular force fields.
  
\subsubsection{Diffusion-Advection Applied to Solvated Proteins}
Note that in this section the variable symbols do not necessarily correspond with those in \textbf{Sec \ref{sec:data}}. Also, this section's argument is heuristic in nature. We proceed by first principles and suppose the solvent moves about the protein by conservation of mass, \textit{ie.} the equation of continuity.
\[
u_t + \Div\left(u\vv\right) = 0
\]
As usual, the material flow will be comprised by diffusion and advection terms. However, to allow that at equilibrium the solvent may not be at uniform density owing to space occupied by the protein, we introduce an effective scaling $\rho(t;x)$ which is identically 1 away from the protein and greater than 1 about the protein. Now we proceed with the usual Fick's Law assertion
\[
u\vv = -\mathbb{K}\nabla\left(u\rho\right) + u\vv_C
\] 
In the simplest conceptual case, the advection velocity $\vv_C$ is the protein's own motion, and the solvent moves by no-slip contact.

By combining these two equations and expanding the gradient, we arrive at the following model which is formally equivalent to (\ref{eq:diff_adv}).
\[
u_t +\Div\left[-\mathbb{K}\rho\nabla u + u(\vv_C - \mathbb{K}\nabla\rho) \right] = 0
\]

\subsection{The Data for the Family of Approximating Problems}\label{sec:data}
We now state more formally the data for the family of approximating diffusion-advection equations over which concentrating capacity is performed. In the process, we allow for the possibility that the diffusivity tensor -- similarly the advection velocity -- is functionally distinct across two regimes: one away from the curve and one at the curve. 

The data for the problems consists of 
\begin{enumerate}
	\item Real numbers $T,\eps_0 > 0$
	\item A smooth and time-varying curve:
	\begin{equation}
	\Gamma:\left[0,T\right]\times \left[-\eps_0,1+\eps_0\right]\rightarrow \Omega\label{eq:Gamma}
	\end{equation}
	\item A positively oriented, and smoothly varying orthonormal triple:
	\begin{align}
	&\vec{t},\vn,\vb:\left[0,T\right]\times \left[-\eps_0,1+\eps_0\right]\rightarrow \rthr
	\label{eq:SO3}\\
	&\begin{pmatrix}
	\vec{t}(t,s) = \partial_s\Gamma/\JG & \vn(t,s) & \vb(t,s)
	\end{pmatrix}\in \mathrm{SO}_{3}(\rr)\nonumber
	\end{align}
	\item A map smoothly parameterizing a moving neighborhood of the curve:
	\begin{align}
	\mathcal{F}:&\left[0,T\right]\times\left[-\eps_0,1+\eps_0\right]\times D_{\eps_0}\rightarrow \left[0,T\right]\times \mathcal{N}_{\eps_0}(t)\label{eq:F}\\
	\mathcal{F}(t,s,\nu,\omega) &=
	\begin{pmatrix}
	t\\
	F(t,s,\nu,\omega) = \Gamma(t,s) + \nu\vn(t,s)+\omega\vb(t,s)\nonumber
	\end{pmatrix}
	\end{align}
	Here $D_{\eps_0}$ is a closed disc in $\rr^2$ of radius $\eps_0$ centered at the origin. $\mathcal{N}_{\eps_0}(t)$ is the image of $F$ at time $t$. We will further denote
	\begin{equation}
	D\mathcal{F} = 
	\begin{pmatrix}
	1 & \vec{0}\\
	\partial_t F & \nabla F \left(=D_{s,\nu,\omega}F\right)
	\end{pmatrix}
	\label{eq:DF}
	\end{equation}
	with respective Jacobians $J_\mathcal{F}$ and $J_F$ as well as the Eulerian velocity field $\vw:\text{Im}\,\mathcal{F}\rightarrow \rr^3$ satisfying
	\begin{equation}
	\vw \circ \mathcal{F} = \partial_t F \label{eq:vw}
	\end{equation} 
	\item Real numbers $\eps\in\left(0,\eps_0\right)$ and $\delta\in\left(0,\eps_0-\eps\right)$
	\item The geometric set and functions
	\begin{align}
	\mathcal{C}_{\eps} &= \left[0,1\right]\times D_{\eps}\label{eq:Ceps}\\
	\chi_{\delta}(t) &= \left(1-t/\delta\right)_+\label{eq:chid}\\
	d_\eps(x) &= \mathrm{dist}(x,\mathcal{C}_\eps) \label{eq:deps}
	\end{align}
	\begin{remark}
		It is evident that $\chi_\delta$ and $d_\eps$ are Lipschitz.
	\end{remark}
	We will also use the notation
	\begin{equation}
	\mathcal{N}_{\eps,\delta}(t) = F(t;\left[d_\eps(x) \leq \delta\right])
	\label{eq:Ned}
	\end{equation}
	\item The decay functions $\zetaed:\left[0,T\right]\times\Omega\rightarrow \rr$ which are (smoothly) zero-extended outside $\mathcal{N}_{\eps_0}(t)$ and inside are determined by the relation
	\begin{equation}
	\zetaed\circ\mathcal{F} = \chi_{\delta}\circ d_{\eps}\label{eq:zetaed}
	\end{equation}
	\begin{remark}
		By the regularity of $\mathcal{F}$,
		$\zetaed$ is Lipschitz because $\chi_\delta$ and $d_\eps$ are.
	\end{remark}
	\item Concentrating capacity coefficients $\aed:\left[0,T\right]\times\bar{\Omega}\rightarrow \rr$ and given by
	\begin{equation}
	\aed(t;x) = 1 + \left(\frac{\eps^2_0}{\eps^2}-1\right)\zetaed(t;x)\label{eq:aed}
	\end{equation}
	\item Diffusivity tensors $\MK_{\eps,\delta}:\left[0,T\right]\times\bar{\Omega}\rightarrow \mathcal{M}_{3\times 3}(\rr)$ given by
	\begin{equation}
	\MK_{\eps,\delta}(t;x) = k_0\mathbb{I}_3 + \left(\frac{\eps^2_0}{\eps^2}\mathcal{K}(t;x) - k_0\mathbb{I}_3\right)\zetaed(t;x)\label{eq:Ked}
	\end{equation}
	$\mathbb{I}_n$ is the n-dimensional identity matrix. We also assume that the principle directions of diffusion for $\mathcal{K}$ are governed by the curve geometry and that
	\begin{itemize}
		\item $\displaystyle 
		\left[\begin{pmatrix}
		\vec{t} & \vn & \vb
		\end{pmatrix}^{\mathrm{t}}\mathcal{K}\begin{pmatrix}
		\vec{t} & \vn & \vb
		\end{pmatrix}\right]\circ \mathcal{F} = 
		\begin{pmatrix}
		k_s(s,\nu,\omega) & \vec{0}\\
		\vec{0} & k_n(s,\nu,\omega)\mathbb{I}_2
		\end{pmatrix}$ 
		\item The functions $k_s,k_n:\left[-\eps_0,1+\eps_0\right]\times D_{\eps_0}\rightarrow \left(0,\infty\right)$ are assumed continuous. 
	\end{itemize}
	\item A real number $\vartheta$ such that
	\begin{equation}
	0 < \vartheta \leq \mathrm{min}\left\{k_0,k_s(s,\nu,\omega),k_n(s,\nu,\omega)\right\}\label{eq:nu}
	\end{equation}
	\item Advection fields $\vved:\left[0,T\right]\times\bar{\Omega}\rightarrow \rthr$ given by
	\begin{equation}
	\vved(t;x) = \vv(t;x) + \left(\frac{\eps^2_0}{\eps^2}\vec{v}_{C}(t;x) - \vec{v}(t;x)\right)\zetaed(t;x)
	\label{eq:vved}
	\end{equation}
	\begin{itemize}
		\item The fields $\vv,\vv_C:\left[0,T\right]\times\bar{\Omega}\rightarrow\rthr$ are assumed continuous.
	\end{itemize}
	\item Initial data $u_0\in\mathcal{C}(\bar{\Omega})$
\end{enumerate}
\subsection{Statement of the Concentrated Capacity Limit}
In the next two sections, we present the equation to be satisfied by the concentrated capacity limit of (\ref{eq:diff_adv},\ref{eq:neu}) and the main theorem providing its small time existence. 
\begin{definition}\label{Def:CC} Defining the Limit \\
	By the concentrated capacity solution of (\ref{eq:diff_adv},\ref{eq:neu}) with data prescribed in \textbf{Sec \ref{sec:data}}, it is meant 
	a solution pair $(\ccu,\ccuC)\in L^2\left(\Oh\right)\times L^2\left(\left[0,1\right]_h\right)$,  whose spatial gradients satisfy $\left(\nabla \ccu, \partial_s \ccuC\right)\in L^2\left(\Oh\right)\times L^2\left(\left[0,1\right]_h\right)$, and satisfying for all smooth test functions $\varphi\in C^{\infty}\left(\bar{\Omega}_h\right)$ such that $\varphi(h;\cdot) \equiv 0$
	\begin{align}
	0 &= \left(- \int_\Omega \varphi(0;\cdot)u_0 dx -\iOh \varphi_t\ccu dxdt  + \iOh \nabla\varphi\cdot k_0\Mbb{I}_3\nabla\ccu dxdt -\iOh\nabla\varphi\ccu\vv dxdt\right)_{\mathrm{Vol}} \label{def:cc}\\
	&+\pi\eps^2_0\left( -\int^1_0 \left(\varphi u_0\right)(0,s,0,0)\JG ds - \int_{\left[0,1\right]_h} \left[\varphi_t - \nabla_{s,\nu,\omega}\varphi\cdot\nabla F^{-1}\partial_t F \right](t,s,0,0)\ccuC \JG dsdt\right. \nonumber\\
	&\;\;\;\;\;\;\;\;\;\;\;\;-\int_{\left[0,1\right]_h} \left[\nabla_{s,\nu,\omega}\varphi\cdot\nabla F^{-1}\vv_C\right](t,s,0,0)\ccuC \JG dsdt \nonumber\\
	&\left.\;\;\;\;\;\;\;\;\;\;\;\;+ \int_{\left[0,1\right]_h} \left[\nabla_{s,\nu,\omega}\varphi\cdot
	\begin{pmatrix}
	k_s/\JG^2 & \vec{0}\\
	\vec{0} & k_n\Mbb{I}_2
	\end{pmatrix}\right](t,s,0,0)
	\begin{pmatrix}
	\partial_s\ccuC\\
	\xi_\nu\\
	\xi_\omega
	\end{pmatrix}
	\JG dsdt\right)_{\mathrm{Line}} \nonumber
	\end{align} 
	The functions $\xi_\nu,\xi_\omega$ are given by
	\begin{align}
	\xi_\nu &= \left[\vec{e}_2\cdot\nabla F^{-1}\left(\vv_C-\partial_tF \right)\right](t,s,0,0)\ccuC/k_n(s,0,0)
	\label{eq:xiwklim}\\
	\xi_\omega &= \left[\vec{e}_3\cdot\nabla F^{-1}\left(\vv_C-\partial_tF \right)\right](t,s,0,0)\ccuC/k_n(s,0,0)\nonumber
	\end{align}
\end{definition}
\begin{remark}
	As will be shown, the form of (\ref{def:cc}) is a natural consequence of the compactness needed to make sense of the limit.
\end{remark}
\begin{remark}
	Observe that the limit is independent of the diffusion coefficient $k_n$. Also, as expected the 1d integrals can be regarded as line integrals taken over $\Gamma$.
\end{remark}
\begin{remark}
	For the volume integrals, $\varphi$ and its gradient are with respect to rectangular coordinates $(t,x,y,z)$ of $\Omega$. For the line integrals, $\varphi$ and its gradient are with respect to local coordinates at the curve $(t,s,\nu,\omega)$.
\end{remark}
\begin{remark}
	Since $\ccuC$ belongs to $L^2\left(\left[0,1\right]_h\right)$ in the variables $(t,s)$, it cannot be differentiated in the variables $\nu,\omega$ orthogonal to the curve. Nevertheless, the limiting process over the family of diffusion-advection equations (\ref{eq:CCFam}) retains a memory of the solution gradient in these directions. It will be shown in a weak sense that  $\left(\xi_\nu,\xi_\omega\right)$ is the orthogonal gradient of $\ccuC$ at the curve.  
\end{remark}
\subsection{The Main Result}
The main argument of this paper will concern the small time existence of concentrated capacity solutions for the diffusion-advection equation (\ref{eq:diff_adv},\ref{eq:neu}).
\begin{theorem}\label{Thm:1}
	For the data prescribed in \textbf{Sec \ref{sec:data}}, there is an $h = h(\text{data}) > 0$ for which a concentrated capacity solution satisfying (\ref{def:cc}) exists.
\end{theorem}

Although the inherent regularity of (\ref{def:cc}) is not sufficient to rigorously justify its strong form, it is instructive to proceed formally with integration by parts and extract a corresponding but formal pointwise PDE from the weak one. We report this here over each subdomain as well as the global domain.

\subsubsection{Formal Pointwise Limit: $\Omega\setminus \Gamma(t)$}
By using an arbitrary test function in (\ref{def:cc}) which vanishes in a space-time neighborhood of the curve and integrating by parts twice, we arrive at a simple diffusion-advection process away from the curve.
\[
\partial_t\ccu - k_0\Delta\ccu + \Div\left(\ccu\vv\right)= 0 \mbox{\;\; in } \Omega\setminus \Gamma(t)
\]

\subsubsection{Formal Pointwise Limit: $\Gamma(t)$}
In (\ref{def:cc}) integrate by parts twice after taking a test function of type $\varphi\zetaed$ with $\varphi$ arbitrary and letting $\delta,\eps\rightarrow 0$. Let $\mathcal{D}_\eps(\Gamma)$ denote the 2d cross section of $\mathcal{N}_{\eps,0}(t)$ centered about $\Gamma(t,s)$. Note that all derivatives below are with respect to the ambient $\left(t,x,y,z\right)$ coordinate system.
\begin{align*}
\partial_t\left(\ccuC\dGt\right) - &\Div\left(\mathcal{K}\nabla\ccuC\dGt\right)+ \Div\left(\ccuC\vv_C\dGt\right)\\
&= \left(\frac{1}{\pi\eps^2_0}\lim_{\eps\rightarrow 0}\int_{\mathcal{D}_\eps(\Gamma)}k_0\Delta\ccu\right)\dGt \mbox{\;\; on } \Gamma(t)
\end{align*}
\begin{remark}
	Because $\Delta\ccu\notin L^1(\Omega_h)$, its limit above need not vanish.
\end{remark}
\subsubsection{Formal Pointwise Limit: Global}
In (\ref{def:cc}) take an arbitrary test function and integrate by parts twice. Again all derivatives below are with respect to the ambient coordinate system.
\begin{align*}
\partial_t\left(\ccu+\pi\eps^2_0\ccuC\dGt \right) - \Div\left(k_0\nabla\ccu + \pi\eps^2_0\mathcal{K}\nabla\ccuC\dGt \right) &+ \Div\left(\ccu\vv+\pi\eps^2_0\ccuC\vv_C\dGt \right) \\
&=0\mbox{\;\; in } \Omega
\end{align*}
\section{Weak Compactness of Approximating Problems}

\subsection{The Family of Approximating Problems}
Consider the following family of diffusion-advection equations with vanishing Neumann data and initial data $u_0$.
\begin{equation}
\partial_t\big(\aed(t;x) u_{\eps,\delta}\big) - \Div\left[\MK_{\eps,\delta}(t;x)\nabla\ued\right] +\Div\left[u\vved(t;x)\right] = 0 \mbox{ in } \Omega_T\label{eq:CCFam}
\end{equation}
The existence of solutions for (\ref{eq:CCFam}) follows from standard theory \cite[Ch. 3]{Ldy68} and properties of the data. Solutions can be taken to satisfy $\ued\in\mathcal{C}\big(\left[0,T\right];L^2(\Omega)\big)\cap L^2\big(\left[0,T\right];W^{1,2}(\Omega)\big)$. 

\begin{remark}
	The placement of the $\left(\eps^2_0/\eps^2\right)\zetaed$ coefficients with respect to the partial derivatives in the first two, left-hand terms of (\ref{eq:CCFam}) is not arbitrary. It has been chosen for mass balance and the requirement that there should exist a compactness estimate like \textbf{Prop \ref{pr:Sobbd}}.
\end{remark}
\subsubsection{\textit{A priori} Estimates}
\begin{proposition}\label{pr:Sobbd}
	For an $h=h(\mathrm{data})$ and $\eps_i\rightarrow 0$, $\exists \kappa_i(\eps_i)(>0)\rightarrow 0$ such that if $\delta_i\leq \kappa_i$ then
	\[
	\esup_{t\leq h}\int_\Omega a_{\eps_i,\delta_i}u^2_{\eps_i,\delta_i}(t;\cdot)dx + \int_{\Omega_h}a_{\eps_i,\delta_i}\left|\nabla u_{\eps_i,\delta_i}\right|^2 dxdt \leq C(\mathrm{data})
	\]
\end{proposition}
\emph{Proof:} \\
For a positive time $h>0$, test by $\ued$ in (\ref{eq:CCFam}) which is justified by Steklov averages \cite[Ch. 2]{DiBenedetto93} and obtain
\begin{align*}
0 &= \int_{\Omega_h}\ued\partial_t\left(\aed\ued\right) dxdt + \int_{\Omega_h}\nabla\ued \MK_{\eps,\delta}\nabla\ued dxdt - \int_{\Omega_h}\nabla\ued\ued\vved dxdt\\
&=\int_{\Omega}\frac{1}{2}\left[\aed\ued^2\right](h;\cdot)dx - \int_{\Omega}\frac{1}{2}\left[\aed u^2_0\right](0;\cdot)dx + \int_{\Omega_h}\frac{1}{2}\aedt\ued^2 dxdt\\
&\;\;+\int_{\Omega_h}\nabla\ued\MK_{\eps,\delta}\nabla\ued dxdt - \int_{\Omega_h}\nabla\ued\ued\vved dxdt
\end{align*}
Rearranging we have
\begin{align}
\int_{\Omega}\frac{1}{2}&\left[\aed\ued^2\right](h;\cdot)dx +\int_{\Omega_h}\nabla\ued\MK_{\eps,\delta}\nabla\ued dxdt \label{eq:1stenest}\\
&\leq \int_{\Omega}\frac{1}{2}\left[\aed u^2_0\right](0;\cdot)dx - \int_{\Omega_h}\frac{1}{2}\aedt\ued^2 dxdt + \int_{\Omega_h}\nabla\ued\ued\vved dxdt \nonumber
\end{align}
From the definitions of $\MK_{\eps,\delta}$ (\ref{eq:Ked}) and $\vartheta$ (\ref{eq:nu}), we may estimate 
\begin{align*}
\nabla\ued^{\mathrm{t}}\MK_{\eps,\delta}\nabla\ued &= k_0\left|\nabla\ued\right|^2\left(1-\zetaed\right) + \frac{\eps^2_0}{\eps^2}\nabla\ued^{\mathrm{t}}\mathcal{K}\nabla\ued\zetaed\\
&\geq \vartheta\left|\nabla\ued\right|^2\left(1-\zetaed\right) + \vartheta\frac{\eps^2_0}{\eps^2}\left|\nabla\ued\right|^2\zetaed\\
&= \vartheta\aed\left|\nabla\ued\right|^2
\end{align*}
Similarly, we may use the definition of $\vved$ (\ref{eq:vved}) to estimate
\begin{align*}
\left|\vved\right| &\leq \left|\vv\right|(1-\zetaed) +\frac{\eps^2_0}{\eps^2}\left|\vv_C\right| \zetaed\\
&\leq \left\|\left\langle\vv,\vv_C\right\rangle\right\|_\infty \aed
\end{align*}
We make these substitutions into (\ref{eq:1stenest}) and arrive at
\begin{align}
\int_{\Omega}\frac{1}{2}&\left[\aed\ued^2\right](h;\cdot)dx +\vartheta\int_{\Omega_h}\aed\left|\nabla\ued\right|^2 dxdt \label{eq:2ndenest}\\
&\leq \int_{\Omega}\frac{1}{2}\left[\aed u^2_0\right](0;\cdot)dx - \int_{\Omega_h}\frac{1}{2}\aedt\ued^2 dxdt + \left\|\left\langle\vv,\vv_C\right\rangle\right\|_\infty \int_{\Omega_h}\aed\left|\nabla\ued\ued\right|dxdt \nonumber
\end{align}
The time derivative of the concentrating capacity coefficients requires some attention
\begin{lemma} The $\aedt$ term\label{lm:aedt}
\[
\aedt\circ\mathcal{F} = \left(\frac{1}{\delta}\left(\frac{\eps^2_0}{\eps^2}-1\right)\nabla d_\eps^{\mathrm{t}}\nabla F^{-1}\partial_t F\right) \chi_{\left[0<d_\eps<\delta\right]}
\]
\end{lemma}
\emph{Proof:} From (\ref{eq:aed}), we see that \[
\aedt = \left(\frac{\eps^2_0}{\eps^2}-1\right)\zetaedt
\]
We may compute $\zetaedt$ from chain-rule and \textbf{Lemma \ref{lm:DFinv}}.
\begin{align*}
\aedt\circ\mathcal{F} &= \left(\frac{\eps^2_0}{\eps^2}-1\right)\zetaedt\circ \mathcal{F} = \left(\frac{\eps^2_0}{\eps^2}-1\right) D\left[\zetaed\circ \mathcal{F}\right]
\begin{pmatrix}
1\\
-\nabla F^{-1}\partial_t F\\
\end{pmatrix}\\
&= \left(\frac{\eps^2_0}{\eps^2}-1\right) D\left[\chi_\delta\circ d_\eps\right]
\begin{pmatrix}
1\\
-\nabla F^{-1}\partial_t F\\
\end{pmatrix}
\end{align*}
We conclude once we observe that
\[
D\left[\chi_\delta\circ d_\eps\right] = 
\begin{pmatrix}
0 & -\frac{1}{\delta}\chi_{\left[0<d_\eps<\delta\right]}\nabla d^{\mathrm{t}}_\eps
\end{pmatrix}
\]
\hfill \bbox\\
Using \textbf{Lemma \ref{lm:aedt}}, a change of variables as in \textbf{Lemma \ref{lm:DF}}, and the coarea formula (\cite[Ch. 2]{Ambrosio00},\cite[Ch. 3]{Evans15}), we may continue
\begin{align*}
\int_{\Omega_h}\aedt\ued^2 dxdt &= \int^h_0\int_{\mathcal{N}_{\eps_0}(t)}\aedt\ued^2 dxdt \\
&= \int^h_0\int_{\left[0<d_\eps<\delta\right]}\frac{1}{\delta}\left(\frac{\eps^2_0}{\eps^2}-1\right)\nabla d_\eps^{\mathrm{t}}\nabla F^{-1}\partial_t F\ued^2 J_F dxdt\\
&= \left(\frac{\eps^2_0}{\eps^2}-1\right)\int^h_0\dashint^\delta_0\int_{\left[d_\eps = \eta\right]}J_F\nabla F^{-1}\partial_t F \ued^2\cdot\frac{\nabla d_\eps}{\left|\nabla d_\eps\right|}d\sigma d\eta dt\\
&= \left(\frac{\eps^2_0}{\eps^2}-1\right)\int^h_0\dashint^\delta_0\int_{\left[d_\eps\leq\eta\right]}\Div_{s,\nu,\omega}\left(J_F\nabla F^{-1}\partial_t F\ued^2\right)dxd\eta dt
\end{align*}
In the last line, we used that $\nabla d_\eps/\left|\nabla d_\eps\right|$ is the exterior normal to $\left[d_\eps\leq\eta\right]$ and applied the divergence theorem. Next, recall that for a vector field $\vec{q}$ defined on $\Omega_T$ coordinates, and which may be shown from a weak formulation,
\begin{equation}\label{eq:divw}
\left(\Div_{x,y,z} \vec{q} \right)\circ F = \frac{1}{J_F}\Div_{s,\nu,\omega}\left[J_F\nabla F^{-1}\left(\vec{q}\circ F\right)\right]
\end{equation}
In particular, we may use this change of variables for the motion field $\vw$, defined by (\ref{eq:vw}). Making this substitution and changing variables back to ambient $x,y,z$ coordinates we have
\[
\int_{\Omega_h}\aedt\ued^2 dxdt = \left(\frac{\eps^2_0}{\eps^2}-1\right)\int^h_0\dashint^\delta_0\int_{\mathcal{N}_{\eps,\eta}(t)}\Div\left(\vw\ued^2\right)dxd\eta dt
\]
We now distribute the divergence, majorize terms, then majorize the integral over $\Nepdl(t)$, and lastly breakup the integral over disjoint domains.
\begin{align}
\left|\int_{\Omega_h}\aedt\ued^2 dxdt\right| & \leq 2\frac{\eps^2_0}{\eps^2}\left\|\left\langle\vw,\nabla\vw\right\rangle\right\|_\infty \int^h_0\int_{\Nepdl(t)}\left(\ued^2 + 2\left|\ued\nabla\ued\right|\right)dx dt \label{eq:3rdenest}\\
&\leq C(\mathrm{data})\left(\int^h_0\int_{\Nepzr(t)}\aed\ued^2 dxdt + \int^h_0\int_{\Nepzr(t)}\aed\left|\ued\nabla\ued\right| dxdt \right.\nonumber\\
&\;\;\;\;\;\;\;\;\;\;\;\;\;\;\;\;\;\;\;\;\;\;\left. +\frac{\eps^2_0}{\eps^2}\int^h_0\int_{\left[\Nepdl\setminus\Nepzr\right](t)} \left(\ued^2+\left|\ued\nabla\ued\right|\right)dxdt\right)\nonumber
\end{align}
We focus now on the terms within $\left[\Nepdl\setminus\Nepzr\right](t)$. Specifically we will look at the $\left|\ued\nabla\ued\right|$ term as the other may be handled similarly. For a still to be determined $\sigma < 2$, apply H\"{o}lder's inequality and \textbf{Lemma \ref{lm:gap}}.
\begin{align*}
\frac{\eps^2_0}{\eps^2}\int^h_0\int_{\left[\Nepdl\setminus\Nepzr\right](t)}&\left|\ued\nabla\ued\right|dxdt\leq \frac{\eps^2_0}{\eps^2} \left\|\ued\left|\nabla\ued\right|\right\|_{\sigma,\Omega_h}\left\|\chi_{\left[0,h\right]\times\left[\Nepdl\setminus\Nepzr\right](t)}\right\|_{\frac{\sigma}{\sigma-1},\Omega_h}\\
&\leq C(\mathrm{data})\frac{\delta^{\frac{\sigma-1}{\sigma}}}{\eps^2} \left\|\ued\left|\nabla\ued\right|\right\|_{\sigma,\Omega_h}\\
&\leq C(\mathrm{data})\frac{\delta^{\frac{\sigma-1}{\sigma}}}{\eps^2} \left(\left\|\left|\nabla\ued\right|^\sigma\right\|_{\frac{2}{\sigma},\Omega_h}\left\|\left|\ued\right|^\sigma\right\|_{\frac{2}{2-\sigma},\Omega_h}\right)^{1/\sigma}\\
&= C(\mathrm{data})\frac{\delta^{\frac{\sigma-1}{\sigma}}}{\eps^2} \left\|\nabla\ued\right\|_{2,\Omega_h}\left\|\ued\right\|_{\frac{2\sigma}{2-\sigma},\Omega_h}
\end{align*}
Recall from the parabolic Sobolev embedding \cite[Ch. 1]{DiBenedetto93}, $\exists \gamma(\Omega,T)$, which is increasing with respect to $T$, so that 
\[
\left\|\ued\right\|_{10/3,\Omega_h}\leq \gamma(\Omega, T)\left(\esup_{t\leq h}\left\|\ued(t;\cdot)\right\|_{2,\Omega_h}+ \left\|\nabla\ued\right\|_{2,\Omega_h}\right)
\]
We now choose $\sigma$ to satisfy $2\sigma/\left(2-\sigma\right)=10/3$ and \textit{a posteriori} verify that $\sigma = 5/4 < 2$ and $\left(\sigma-1\right)/\sigma = 1/5$. For this choice we next apply Young's inequality and majorize with $1\leq\aed$.
\begin{align}
\frac{\eps^2_0}{\eps^2}\int^h_0\int_{\left[\Nepdl\setminus\Nepzr\right](t)}&\left|\ued\nabla\ued\right|dxdt \label{eq:gapokudu}\\
&\leq C(\mathrm{data})\frac{\delta^{1/5}}{\eps^2}\left(\esup_{t\leq h}\int_{\Omega}\aed\ued^2\left(t;\cdot\right)dx + \int_{\Omega_h}\aed\left|\nabla\ued\right|^2 dxdt \right)\nonumber
\end{align}
By writing $\ued^2 = \left|\ued\ued\right|$ and majorizing $\left\|\ued\right\|^2_{2,\Omega_h} \leq h\left(\esup_{t\leq h}\int_{\Omega}\ued^2(t;\cdot)dx\right)$, we may proceed similarly as when it was $\left|\ued\nabla\ued\right|$. 
\begin{align}
\frac{\eps^2_0}{\eps^2}\int^h_0\int_{\left[\Nepdl\setminus\Nepzr\right](t)}&\ued^2 dxdt \label{eq:gapokuu} \leq C(\mathrm{data})\frac{\delta^{\frac{\sigma-1}{\sigma}}}{\eps^2} \left\|\ued\right\|_{2,\Omega_h}\left\|\ued\right\|_{\frac{2\sigma}{2-\sigma},\Omega_h} \\
&\leq C(\mathrm{data})\frac{\delta^{1/5}}{\eps^2}\left(\esup_{t\leq h}\int_{\Omega}\aed\ued^2\left(t;\cdot\right)dx + \int_{\Omega_h}\aed\left|\nabla\ued\right|^2 dxdt \right)\nonumber
\end{align}
Now we return to the energy terms of (\ref{eq:3rdenest}) in $\Nepzr(t)$. These are simpler to estimate. Let $\alpha > 0$ and use Young's inequality.
\begin{align}
\int^h_0\int_{\Nepzr(t)}\aed\ued^2 dxdt &\leq h\left(\esup_{t\leq h}\int_\Omega\aed\ued^2(t;\cdot)dx \right)\label{eq:ngapuu}\\
\int^h_0\int_{\Nepzr(t)}\aed\left|\ued\nabla\ued\right|dxdt &\leq \frac{\alpha}{2}\int_{\Omega_h}\aed\left|\nabla\ued\right|^2dxdt + \frac{h}{2\alpha}\left(\esup_{t\leq h}\int_{\Omega}\aed\left|\ued\right|^2\left(t;\cdot\right)dx \right)\label{eq:ngapudu}
\end{align}
We finally combine (\ref{eq:gapokudu},\ref{eq:gapokuu},\ref{eq:ngapuu},\ref{eq:ngapudu}) into (\ref{eq:3rdenest}).
\begin{align}
\left|\int_{\Omega_h}\aedt\ued^2 dxdt\right| \leq C(\mathrm{data})&\left[ \left(h + \frac{h}{2\alpha}+\frac{\delta^{1/5}}{\eps^2}\right)\left(\esup_{t\leq h}\int_\Omega\aed\left|\ued\right|^2\left(t;\cdot\right)dx\right)\right.\label{eq:aedtenest}\\
&\left. \left(\frac{\alpha}{2}+\frac{\delta^{1/5}}{\eps^2}\right)\int_{\Omega_h}\aed\left|\nabla\ued\right|^2dxdt\right]\nonumber 
\end{align}
Having estimated the $\aedt\ued^2$ energy term, we next look at the advection term of (\ref{eq:2ndenest}). This is estimated by a simple Young's inequality.
\begin{equation}\label{eq:advenest}
\int_{\Omega_h}\aed\left|\ued\nabla\ued\right|dxdt \leq \left(\frac{\alpha}{2}\int_{\Omega_h}\aed\left|\nabla\ued\right|^2 dxdt + \frac{h}{2\alpha}\esup_{t\leq h}\int_{\Omega}\left[\aed\ued^2\right](t;\cdot)dx\right)
\end{equation}
We are now ready to return to (\ref{eq:2ndenest}) and begin concluding the uniform Sobolev estimate by applying (\ref{eq:aedtenest},\ref{eq:advenest}).
\begin{align}
\int_{\Omega}\frac{1}{2}&\left[\aed\ued^2\right](h;\cdot)dx +\vartheta\int_{\Omega_h}\aed\left|\nabla\ued\right|^2 dxdt \label{eq:4thenest}\\
&\leq C(\mathrm{data})\left[ \left(h + \frac{h}{2\alpha}+\frac{\delta^{1/5}}{\eps^2}\right)\left(\esup_{t\leq h}\int_\Omega\aed\left|\ued\right|^2\left(t;\cdot\right)dx\right)\right.\nonumber\\
&\;\;\;\;\;\;\;\;\;\;\;\;\;\;\;\;\;\;\;\left. \left(\frac{\alpha}{2}+\frac{\delta^{1/5}}{\eps^2}\right)\int_{\Omega_h}\aed\left|\nabla\ued\right|^2dxdt\right] + \int_\Omega \frac{1}{2}\left[\aed u^2_0\right]\left(0;\cdot\right)dx\nonumber
\end{align}
As the right side is an increasing function of $h$, we may add an $\esup_{t\leq h}$ to the time trace with argument $(t;\cdot)$ on the left-hand side at the small expense of increasing $C(\mathrm{data})$. It is now clear that we may pick $\alpha = \alpha(\mathrm{data})$, $h = h(\alpha,\mathrm{data})$, and $\delta \leq C(\alpha,h,\mathrm{data})\eps^{11}$, so that for some new $\tilde{C}(\mathrm{data})$,
\[
\esup_{t\leq h}\int_\Omega\left[\aed\ued^2\right]\left(t;\cdot\right)dx + \int_{\Omega_h}\aed\left|\nabla\ued\right|^2 dxdt \leq \tilde{C}(\mathrm{data})\int_\Omega\left[\aed u^2_0\right]\left(0;\cdot\right)dx
\]
By \textbf{Lemma \ref{lm:aedwklim}}, we see that once a sequence $\eps_i\rightarrow 0$ is chosen that there is a $\kappa_i$ so that if $\delta_i\leq \kappa_i$, the right hand integral has a limit and so the whole sequence is bounded. In particular, the $\delta = O(\eps^{11})$ choice is sufficient.\hfill \bbox
\subsection{Extracting Weakly Convergent Subsequences}
From \textbf{Prop \ref{pr:Sobbd}}, we may conclude the existence of weak limits in both the volume and at the curve. 
\begin{remark}
	By abuse of notation, we will continue to index the convergent subsequences with the same single index notation for which the energy estimates were shown.
\end{remark}
\begin{proposition}\label{pr:wkconv}
	For $\eps_i\rightarrow 0$ and $\delta_i < \kappa_i$ as in \textbf{Prop \ref{pr:Sobbd}}, we may find a further subsequence (redundantly notated by $i$) such that weakly
	\begin{align*}
	\ueidi \rightharpoonup \ccu &\mbox{ and } \nabla\ueidi \rightharpoonup \nabla\ccu \mbox{ in } L^2(\Omega_h)\\
	\dashint_{D_{\eps_i}}\ueidi(t,s,\nu,\omega)d\nu d\omega \rightharpoonup \ccuC &\mbox{ and } \dashint_{D_{\eps_i}} \nabla_{s,\nu,\omega}\ueidi d\nu d\omega \rightharpoonup 
	\begin{pmatrix}
	\partial_s \ccuC \\
	\xi_\nu \\
	\xi_\omega
	\end{pmatrix}
	\mbox{ in } L^2\left(\left[0,1\right]_h\right)
	\end{align*}
\end{proposition}
\begin{remark}
	The weak limits $\xi_v,\xi_w\in L^2\left(\left[0,1\right]_h\right)$ with initially unknown relationship to $\ccuC$ 
	emerge because the curve is intrinsically one dimensional, two fewer dimensions than the ambient 3d space. Their final determination is deferred to \textbf{Sec. \ref{S:4}}. See also \textbf{Remark \ref{rmk:wgradlim}} and \textbf{Remark \ref{rmk:wgradlimG}} below.
\end{remark}
\subsubsection{Convergence Over the Volume}
Since $\aed\geq 1$, we have the basic estimate
\[
\esup_{t\leq h}\int_\Omega u^2_{\eps_i,\delta_i}\left(t;\cdot\right)dx + \int_{\Omega_h}\left|\nabla u_{\eps_i,\delta_i}\right|^2 dxdt \leq C(\mathrm{data})
\]
Accordingly, the $\ueidi$ and $\left|\nabla\ueidi\right|$ are uniformly bounded in $L^2(\Omega_h)$ and we may extract weakly convergent subsequences. In particular
\[
\ueidi \rightharpoonup \ccu \mbox{ and } \nabla\ueidi \rightharpoonup \nabla\ccu \mbox{ in } L^2(\Omega_h)
\]
\begin{remark}\label{rmk:wgradlim}
	Recall that the weak limit of the gradients tends to the gradient of the weak limit by a simple integration by parts argument. If $\ueidi \rightharpoonup \ccu $ and $\partial_x\ueidi \rightharpoonup \phi$, then for an arbitrary, compactly supported and smooth test function $\varphi$ it holds
	\[
	\iOh \varphi \phi dxdt = \lim_i \iOh \varphi \partial_x\ueidi dxdt = \lim_i -\iOh \partial_x\varphi \ueidi dxdt = -\iOh \partial_x\varphi \ccu dxdt
	\]
	It follows $\phi = \partial_x\ccu$. Observe this argument is not possible unless the domain $\Omega$ contains the direction for which integration by parts is to be performed. 
\end{remark}
\subsubsection{Convergence Over the Curve}
Restrict attention to the subdomain $\Nepizr(t)$ where ${\aeidi}_{\left.\right| \Nepizr(t)} \equiv \frac{\eps^2_0}{\eps^2_i}$. Here the estimates of \textbf{Prop \ref{pr:Sobbd}} yield 
\[
\esup_{t\leq h}\frac{\eps^2_0}{\eps^2_i}\int_{\Nepizr(t)}\ueidi^2\left(t;\cdot\right)dx + \frac{\eps^2_0}{\eps^2_i}\int^h_0\int_{\Nepizr(t)}\left|\nabla\ueidi\right|^2 dxdt \leq C(\mathrm{data})
\]
Now perform a change of variables to a common domain and conclude with elementary manipulations
\begin{align}
\esup_{t\leq h}\int^1_{s=0}\dashint_{D_{\eps_i}}\ueidi^2(t,s,\nu,\omega)J_F dx &+ \int^h_0\int^1_{s=0}\dashint_{D_{\eps_i}}\nabla_{s,\nu,\omega}\ueidi\cdot\left(\nabla F^{\mathrm{t}}\nabla F\right)^{-1}\nabla_{s,\nu,\omega}\ueidi J_F dxdt \label{eq:1stGcmpt}\\
&\leq C(\mathrm{data}) \nonumber
\end{align}
\begin{lemma}\label{lm:FFcoerc}
	$\exists\beta(\mathrm{data}) > 0$ such that
	\[
	\beta\left|\vec{q}\right|^2 \leq \vec{q}\cdot \left[\left(\nabla F^{\mathrm{t}}\nabla F\right)^{-1}J_F\right]\vec{q}
	\]
\end{lemma}
\emph{Proof:} It is sufficient to show for a unit vector $\vec{q}$. From \textbf{Lemma \ref{lm:DFinv}}, we see that the matrix is symmetric and so must have an orthonormal basis of eigenvectors. Because the matrix is invertible, these eigenvalues are nonzero. Since further, the quadratic form may be written as $\vec{p}\cdot \vec{p} J_F \geq 0$ for $p=\left(\nabla F^{-1}\right)^{\mathrm{t}}\vec{q}$, the eigenvalues are nonnegative. Altogether the eigenvalues are all strictly positive. 

With the variational characterization of the minimum eigenvalue, $\lambda_{\min}$, for a symmetric matrix $A$
\[
\lambda_{\min} = \min_{\left|u\right|=1} u^{\mathrm{t}}Au
\] 
we may consider the map
\[
(\vec{q},x)\rightarrow \vec{q}^\tr\left(\left[\nabla F^\tr \nabla F\right](x)\right)^{-1}\vec{q}:\left[\left|\vec{q}\right|=1\right]\times\left[-\eps_0,1+\eps_0\right]\times D_{\eps_0} \rightarrow \left(0,\infty\right)
\]
The map is continuous, its domain is compact, and is strictly positive due to the considerations of the last paragraph. Hence, its global minimum is also strictly positive, being attained. Since $J_F(x)\geq \gamma > 0$ for some $\gamma$ by assumptions made on the data, the result is concluded.
\hfill \bbox

Now minorize $J_F$, combine \textbf{Lemma \ref{lm:FFcoerc}} with (\ref{eq:1stGcmpt}), and apply Jensen's inequality.
\begin{align}
\esup_{t\leq h}\int^1_{s=0}\left(\dashint_{D_{\eps_i}} \ueidi d\nu d\omega\right)^2ds &+ \int^h_0\int^1_{s=0}\left(\dashint_{D_{\eps_i}}\left|\nabla_{s,\nu,\omega}\ueidi\right|d\nu d\omega\right)^2dsdt \label{eq:2ndGcmpt}\\
&\leq C(\mathrm{data})\nonumber
\end{align}
It now follows that the $\dashint_{D_{\eps_i}} \ueidi d\nu d\omega$ and $\dashint_{D_{\eps_i}}\nabla_{s,\nu,\omega}\ueidi d\nu d\omega$ are uniformly bounded in $L^2\left(\left[0,1\right]_h\right)$, and we may extract weakly convergent subsequences
\[
\dashint_{D_{\eps_i}}\ueidi d\nu d\omega \rightharpoonup \ccuC \mbox{ and } \dashint_{D_{\eps_i}}\nabla_{s,\nu,\omega}\ueidi d\nu d\omega \rightharpoonup 
\begin{pmatrix}
\partial_s \ccuC \\
\xi_\nu\\
\xi_\omega
\end{pmatrix}
\mbox{ in } L^2(\left[0,1\right]_h)
\]
\begin{remark}\label{rmk:wgradlimG}
	Observe that the domain $\left[0,1\right]_h$ in variables $(t,s)$ does contain the unit $s$-direction. Consistent with \textbf{Remark \ref{rmk:wgradlim}}, the corresponding weak limit may be identified with $\partial_s\ccu$. However, it does not contain the unit directions for $\nu$ and $\omega$. For these cases compactness only implies the membership of $\xi_\nu,\xi_\omega$ in $L^2\left(\left[0,1\right]_h\right)$.
\end{remark}
\section{Determining the Concentrated Capacity Limit}
We now proceed with computing the equation satisfied by the concentrated capacity limit. First, we will prove \textbf{Theorem \ref{Thm:1}} without yet determining the unknown weak limits $\xi_\nu,\xi_\omega$. 
\begin{remark}
	To not overburden notation, hereafter we will drop the subscript $i$ from solutions but understand we work with a sequence given by \textbf{Prop \ref{pr:wkconv}}. 
\end{remark}
\begin{proposition}
	The solution pair $\left(\ccu,\ccuC\right)$ provided by \textbf{Prop \ref{pr:wkconv}} satisfies the integral, weak formulation (\ref{def:cc}).
\end{proposition}
\emph{Proof:} For a smooth test function $\varphi\in \mathcal{C}^\infty(\bar{\Omega}_h)$ and vanishing at time $t = h$, pass to the weak formulation of (\ref{eq:CCFam}). 
\begin{align*}
0 &= \int_{\Omega_h}\varphi\bigg(\partial_t\big(\aed u_{\eps,\delta}) - \Div \left[\MK_{\eps,\delta}\nabla\ued\right] + \Div\left[\ued\vved\right]\bigg)dxdt \\
&= -\int_{\Omega_h}\varphi_t \aed\ued dxdt - \int_\Omega \varphi\aed(0;\cdot)u_0 dx\\
&\;\;\;\;+\iOh \nabla\varphi\cdot\MK_{\eps,\delta}\nabla\ued dxdt - \iOh \nabla\varphi\cdot\ued\vved dxdt
\end{align*}
Next we group these terms by those in the volume and those about the curve.
\begin{align}
0 = &\left(-\iOh \varphi_t \left(1-\zetaed\right)\ued dxdt - \int_\Omega \left[\varphi\left(1-\zetaed\right)\right]\left(0;\cdot\right) u_0 dx\right. \label{eq:1stwklim}\\
&\;\;+\left. \iOh \nabla\varphi\cdot k_0 \Mbb{I}_3 \left(1-\zetaed\right)\nabla\ued dxdt - \iOh \nabla\varphi\ued\vv\left(1-\zetaed\right)dxdt\right)_{\mathrm{Vol}} \nonumber\\
+&\left(-\iOh\varphi_t \left(\frac{\eps^2_0}{\eps^2}\zetaed\right)\ued dxdt - \int_\Omega \left[\varphi \left(\frac{\eps^2_0}{\eps^2}\zetaed\right)\right]\left(0;\cdot\right)u_0 dx \right. \nonumber \\
&\;\,+ \left. \iOh\nabla\varphi\cdot \left(\frac{\eps^2_0}{\eps^2}\mathcal{K}\zetaed\right)\nabla\ued dxdt - \iOh \nabla\varphi \ued\left(\frac{\eps^2_0}{\eps^2}\vv_C\zetaed\right)dxdt \right)_{\mathrm{Line}} \nonumber
\end{align}
We will now examine the terms in subsets.
\subsection{Limits in the Volume}
Because $1-\zetaed\rightarrow 1$ strongly in $L^2(\Omega_h)$, we may apply this in tandem with the volume weak limit of \textbf{Prop \ref{pr:wkconv}}. The limit becomes
\begin{equation}\label{eq:wklimvol}
\left(-\iOh \varphi_t\ccu dxdt - \int_\Omega \varphi(0;\cdot)u_0 dx + \iOh \nabla\varphi\cdot k_0\Mbb{I}_3\nabla\ccu dxdt -\iOh\nabla\varphi\ccu\vv dxdt\right)_{\mathrm{Vol}}
\end{equation}
\subsection{Limits on the Curve}
\subsubsection{The Initial Time Term}
As a direct application of \textbf{Lemma \ref{lm:aedwklim}}, we may determine the limit of the initial time data on the curve.
\begin{equation}\label{eq:wklimlineu0}
-\int_\Omega\left[\varphi\left(\frac{\eps^2_0}{\eps^2}\zetaed\right)\right]\left(0;\cdot\right)u_0 dx \rightarrow -\pi\eps^2_0 \int^1_0 \left(\varphi u_0\right)(0,s,0,0)\JG ds
\end{equation}
\begin{remark}
	\textbf{Lemma \ref{lm:aedwklim}} was stated for the $\aed$ coefficients, but its proof showed convergence of this term as well.
\end{remark}
\subsubsection{Controlling the Remaining $\left[\Nepdl\setminus\Nepzr\right](t)$ Terms}
For the remaining integrals in $\left(\ldots\right)_{\text{Line}}$, the contribution due to $\left[\Nepdl\setminus\Nepzr\right](t)$ is vanishing. More precisely, we have the following. 
\begin{lemma}\label{lm:gapwklimest}
\begin{align*}
&\left|\int^h_0\int_{\left[\Nepdl\setminus\Nepzr\right](t)}\varphi_t\left(\frac{\eps^2_0}{\eps^2}\zetaed\right)\ued dxdt\right|,\left|\int^h_0\int_{\left[\Nepdl\setminus\Nepzr\right](t)}\nabla\varphi\cdot\left(\frac{\eps^2_0}{\eps^2}\mathcal{K}\zetaed\right)\nabla\ued dxdt\right|, \\
&\;\;\;\;\; \left|\int^h_0\int_{\left[\Nepdl\setminus\Nepzr\right](t)} \nabla\varphi \ued \left(\frac{\eps^2_0}{\eps^2}\vv_C\zetaed\right)dxdt\right| \leq C(\mathrm{data})\left\|\left\langle\varphi,\varphi_t,\nabla\varphi\right\rangle\right\|_\infty \left(\frac{\delta^{1/2}}{\eps^2}\right)
\end{align*}
\end{lemma}
\emph{Proof:} We handle the first term, the others being similar. 
\begin{align*}
\left|\int^h_0\right.&\int_{\left[\Nepdl\setminus\Nepzr\right](t)}\varphi_t\left.\left(\frac{\eps^2_0}{\eps^2}\zetaed\right)\ued dxdt\right| \leq \frac{\eps^2_0}{\eps^2}\left\|\varphi_t\right\|_{\infty}\int^h_0\int_{\left[\Nepdl\setminus\Nepzr\right](t)}\left|\ued\right|dxdt \\
&\leq \frac{\eps^2_0}{\eps^2}\left\|\varphi_t\right\|_\infty\left\|\ued\right\|_{2,\Omega_h}\left\|\chi_{\left[0,h\right]\times\left[\Nepdl\setminus \Nepzr\right](t)}\right\|_{2,\Omega_h}
\end{align*}
We may use that $1\leq\aed$ and \textbf{Prop \ref{pr:Sobbd}} to bound the energy norm of $\ued$ by a constant depending on the data. Further, by \textbf{Lemma \ref{lm:gap}}, the $L^2$ space time norm of the gap must be bounded like $C(\mathrm{data})\delta^{1/2}$. \hfill \bbox
\subsubsection{The $\varphi_t$ and Advection Terms} Use \textbf{Lemma \ref{lm:gapwklimest}} and a change of variables to obtain
\begin{align*}
-&\iOh \varphi_t\left(\frac{\eps^2_0}{\eps^2}\zetaed\right)\ued dxdt - \iOh \nabla\varphi\ued\left(\frac{\eps^2_0}{\eps^2}\vv_C\zetaed\right)dxdt\\
&= -\int^h_0\int_{\Nepzr(t)}\varphi_t\left(\frac{\eps^2_0}{\eps^2}\right)\ued dxdt - \int^h_0\int_{\Nepzr(t)}\nabla\varphi\ued\left(\frac{\eps^2_0}{\eps^2}\vv_C\right)dxdt + O\left(\frac{\delta^{1/2}}{\eps^2}\right)\\
&= -\pi\eps^2_0\int^h_0\int^1_0\dashint_{D_\eps}\varphi_t\ued J_F d\nu d\omega dsdt - \pi\eps^2_0\int^h_0\int^1_0\dashint_{D_\eps}\nabla\varphi\ued\vv_C J_F d\nu d\omega dsdt + O\left(\frac{\delta^{1/2}}{\eps^2}\right)
\end{align*}
Note that the $\varphi_t,\nabla\varphi$ here are derivatives of $\varphi(t,x,y,z)$ and then composed with $\mathcal{F}$ but not yet the derivatives of $\varphi(t,s,\nu,\omega)$. We now apply chain-rule with the help of \textbf{Lemma \ref{lm:DFinv}} to revert to the coordinate system at the curve.
\begin{align*}
-\iOh \varphi_t&\left(\frac{\eps^2_0}{\eps^2}\zetaed\right)\ued dxdt - \iOh \nabla\varphi\ued\left(\frac{\eps^2_0}{\eps^2}\vv_C\zetaed\right)dxdt\\
&= -\pi\eps^2_0\int^h_0\int^1_0\dashint_{D_\eps}\left(\varphi_t - \nabla_{s,\nu,\omega}\varphi\cdot\nabla F^{-1}\partial_t F \right)\ued J_F d\nu d\omega dsdt \\
&\;\;\;\;-\pi\eps^2_0\int^h_0\int^1_0\dashint_{D_\eps}\nabla_{s,\nu,\omega}\varphi\cdot \nabla F^{-1}\ued\vv_C J_F d\nu d\omega dsdt + O\left(\frac{\delta^{1/2}}{\eps^2}\right)
\end{align*}
From the uniform continuity over $\text{Dom}\; \mathcal{F}$ of all terms except for $\ued$ whose integral average converges weakly by \textbf{Prop \ref{pr:wkconv}} and recalling $\delta = O(\eps^{11})$, we have
\begin{align}
-\iOh& \varphi_t\left(\frac{\eps^2_0}{\eps^2}\zetaed\right)\ued dxdt - \iOh \nabla\varphi\ued\left(\frac{\eps^2_0}{\eps^2}\vv_C\zetaed\right)dxdt \label{eq:wklimlineadv}\\
&\rightarrow -\pi\eps^2_0 \left(\int_{\left[0,1\right]_h} \left[\varphi_t - \nabla_{s,\nu,\omega}\varphi\cdot\nabla F^{-1}\partial_t F \right](t,s,0,0)\ccuC\JG dsdt \right.\nonumber\\
&\;\;\;\;\;\;\;\;\;\;\;\;\;\;\;\;\left. + \int_{\left[0,1\right]_h} \left[\nabla_{s,\nu,\omega}\varphi\cdot\nabla F^{-1}\vv_C\right](t,s,0,0)\ccuC \JG dsdt\right)
\nonumber
\end{align}
\subsubsection{The Diffusion Term}
By \textbf{Lemma \ref{lm:gapwklimest}} and a change of variables we have
\begin{align*}
\iOh& \nabla\varphi\cdot\left(\frac{\eps^2_0}{\eps^2}\mathcal{K}\zetaed\right)\nabla\ued dxdt = \int^h_0\int_{\Nepzr(t)}\nabla\varphi\cdot\left(\frac{\eps^2_0}{\eps^2}\mathcal{K}\right)\nabla\ued dxdt + O\left(\frac{\delta^{1/2}}{\eps^2}\right)\\
&= \pi\eps^2_0\int^h_0\int^1_0\dashint_{D_\eps}\nabla_{s,\nu,\omega}\varphi\cdot\left[\nabla F^{-1}\mathcal{K}\left(\nabla F^{-1}\right)^\tr\right]\nabla_{s,\nu,\omega}\ued J_F d\nu d\omega dsdt + O\left(\frac{\delta^{1/2}}{\eps^2}\right)
\end{align*}
Appealing to the uniform continuity over $\text{Dom}\; \mathcal{F}$ of all terms save $\nabla_{s,\nu,\omega}\ued$ whose integral average converges weakly by \textbf{Prop \ref{pr:wkconv}} and recalling that $\delta = O(\eps^{11})$, we have
\begin{align*}
\iOh& \nabla\varphi\cdot\left(\frac{\eps^2_0}{\eps^2}\mathcal{K}\zetaed\right)\nabla\ued dxdt \\
&\rightarrow \pi\eps^2_0 \int^h_0\int^1_0 \left[\nabla_{s,\nu,\omega}\varphi\cdot\nabla F^{-1}\mathcal{K}\left(\nabla F^{-1}\right)^\tr\right](t,s,0,0)
\begin{pmatrix}
\partial_s\ccuC\\
\xi_\nu\\
\xi_\omega
\end{pmatrix}
\JG dsdt
\end{align*}
From \textbf{Lemma \ref{lm:DFinv}} and the structure conditions of $\mathcal{K}$ (\ref{eq:Ked}) we have
\[
\left[\nabla F^{-1}\right](t,s,0,0) = 
\begin{pmatrix}
\JG^{-2} & 0 & 0\\
0 & 1 & 0\\
0 & 0 & 1
\end{pmatrix}
\begin{pmatrix}
\partial_s\Gamma^\tr \\
\vn^\tr \\
\vb^\tr
\end{pmatrix}
\]
from whence it follows
\[
\left[\nabla F^{-1}\mathcal{K}\left(\nabla F^{-1}\right)^\tr\right](t,s,0,0) = \begin{pmatrix}
k_s(s,0,0)/\JG^2 & 0 & 0\\
0 & k_n(s,0,0) & 0\\
0 & 0 & k_n(s,0,0)
\end{pmatrix}
\]
Combining we have
\begin{align}
\iOh& \nabla\varphi\cdot\left(\frac{\eps^2_0}{\eps^2}\mathcal{K}\zetaed\right)\nabla\ued dxdt \label{eq:wklimdiff}\\
&\rightarrow \pi\eps^2_0 \int^h_0\int^1_0 \left[\nabla_{s,\nu,\omega}\varphi\cdot
\begin{pmatrix}
k_s/\JG^2 & \vec{0}\\
\vec{0} & k_n\Mbb{I}_2
\end{pmatrix}\right](t,s,0,0)
\begin{pmatrix}
\partial_s\ccuC\\
\xi_\nu\\
\xi_\omega
\end{pmatrix}
\JG dsdt\nonumber
\end{align}
\subsection{Concluding the Limit}
We now conclude using (\ref{eq:wklimvol},\ref{eq:wklimlineu0},\ref{eq:wklimlineadv},\ref{eq:wklimdiff}) that (\ref{eq:1stwklim}) limits to
\begin{align}
0 &= \left(- \int_\Omega \varphi(0;\cdot)u_0 dx -\iOh \varphi_t\ccu dxdt  + \iOh \nabla\varphi\cdot k_0\Mbb{I}_3\nabla\ccu dxdt -\iOh\nabla\varphi\ccu\vv dxdt\right)_{\mathrm{Vol}} \label{eq:wklim}\\
&+\pi\eps^2_0\left( -\int^1_0 \left(\varphi u_0\right)(0,s,0,0)\JG ds - \int_{\left[0,1\right]_h} \left[\varphi_t - \nabla_{s,\nu,\omega}\varphi\cdot\nabla F^{-1}\partial_t F \right](t,s,0,0)\ccuC \JG dsdt\right. \nonumber\\
&\;\;\;\;\;\;\;\;\;\;\;\;-\int_{\left[0,1\right]_h} \left[\nabla_{s,\nu,\omega}\varphi\cdot\nabla F^{-1}\vv_C\right](t,s,0,0)\ccuC \JG dsdt \nonumber\\
&\left.\;\;\;\;\;\;\;\;\;\;\;\;+ \int_{\left[0,1\right]_h} \left[\nabla_{s,\nu,\omega}\varphi\cdot
\begin{pmatrix}
k_s/\JG^2 & \vec{0}\\
\vec{0} & k_n\Mbb{I}_2
\end{pmatrix}\right](t,s,0,0)
\begin{pmatrix}
\partial_s\ccuC\\
\xi_\nu\\
\xi_\omega
\end{pmatrix}
\JG dsdt\right)_{\mathrm{Line}} \nonumber
\end{align}
\hfill \bbox
\section{Determining the Weak Limits $\xi_v,\xi_w$}\label{S:4}
To conclude the expressions for $\xi_v,\xi_w$ we will need to use test functions of the form 
\begin{equation}\label{eq:testvphi}
\varphi_{\eps,\delta} = \psi+\left(\psi_C-\psi \right)\zetaed
\end{equation}
for $\psi,\psi_C\in\mathcal{C}^\infty(\bar{\Omega}_h)$ and vanish at time $t=h$ in (\ref{eq:wklim}).  {The use of $\zetaed$ is justified by a simple mollification argument.} Because $\varphi_{\eps,\delta}\equiv \psi_C$ in a neighborhood of the curve, we may directly substitute $\psi_C$ for $\varphi$ in $\left(\ldots\right)_{\mathrm{Line}}$ of (\ref{eq:wklim}). The corresponding $\left(\ldots\right)_{\mathrm{Vol}}$ terms need consideration.  Note that $\varphi_\delta := \varphi_{\delta,\delta}$.
\begin{lemma}\label{lm:testvphi}
	For the choice of test function (\ref{eq:testvphi}) with $\eps=\delta$ and two smooth test functions $\psi,\psi_C$ which uniformly satisfy for all times $t$, $\left\|\psi-\psi_C\right\|_{\infty,\mathcal{N_{\delta,\delta}}(t)} = O(\delta)$, it holds that as $\delta\rightarrow 0$ the volume terms in (\ref{eq:wklim}) have the limit
	\begin{align*}
	-\iOh &\partial_t\varphi_{\delta}\ccu dxdt - \int_\Omega \varphi_{\delta}(0;\cdot)u_0 dx + \iOh \nabla\varphi_\delta\cdot k_0\Mbb{I}_3\nabla\ccu dxdt -\iOh\nabla\varphi_\delta\ccu\vv dxdt\\
	&\rightarrow -\iOh \partial_t\psi\ccu dxdt - \int_\Omega \psi(0;\cdot)u_0 dx + \iOh \nabla\psi\cdot k_0\Mbb{I}_3\nabla\ccu dxdt -\iOh\nabla\psi\ccu\vv dxdt
	\end{align*}
\end{lemma}
\emph{Proof:} Examine first the diffusion term.
\begin{align*}
\iOh \nabla\varphi_\delta\cdot k_0\Mbb{I}_3\nabla\ccu dxdt &= \iOh \left[\nabla\psi + \nabla\left(\psi_C-\psi\right)\zetadd +  \left(\psi_C-\psi\right)\nabla\zetadd\right]\cdot k_0\Mbb{I}_3\nabla\ccu dxdt \\
&= \iOh \left[\nabla\psi + \nabla\left(\psi_C-\psi\right)\zetadd\right]\cdot k_0\Mbb{I}_3\nabla\ccu dxdt \\
&\;\; - \frac{1}{\delta}\int^h_0\int_{\left[0<d_\delta<\delta\right]}\left[\left(\psi_C-\psi\right)\right](t,s,\nu,\omega)\nabla d_\eps\cdot k_0\left(\nabla F^\tr\nabla F\right)^{-1}\nabla_{s,\nu,\omega}\ccu J_Fdxdt
\end{align*}
By assumption $\left|\psi_C-\psi\right| = O(\delta)$ over $\left[d_\delta < \delta\right] $, and so the $\delta^{-1}$ singularity is balanced. The remaining integral vanishes in the limit from the absolute continuity of $\nabla\ccu$ and the sup bounds on the remaining terms. Since also $\zeta_{\delta,\delta}\rightarrow 0$ a.e., this concludes the first term. 

The time derivative and advection terms follow similarly using the absolute continuity of $\ccu$. For the former we have below. 
\begin{align*}
-\iOh \partial_t\varphi_{\delta}\ccu dxdt &= -\iOh \left[\partial_t\psi +\left(\partial_t\psi_{C}-\partial_t\psi\right)\zetadd + \left(\psi_C-\psi\right)\zetaddt\right]\ccu dxdt \\
&= -\iOh \left[\partial_t\psi +\left(\partial_t\psi_{C}-\partial_t\psi\right)\zetadd\right]\ccu dxdt \\
&\;\;\;\, -\frac{1}{\delta} \int^h_0\int_{\left[0<d_\eps<\delta\right]}\left[\left(\psi_C-\psi\right)\ccu\right](t,s,\nu,\omega)\nabla d_\eps\cdot \nabla F^{-1}\partial_t F J_F dx dt
\end{align*}
For the latter we have below.
\begin{align*}
-\iOh \nabla\varphi_\delta\ccu\vv dxdt &= -\iOh \left[\nabla\psi + \nabla\left(\psi_C-\psi\right)\zetadd + \left(\psi_C-\psi\right)\nabla\zetadd\right]\ccu\vv dxdt \\
&= -\iOh\left[ \nabla\psi + \nabla\left(\psi_C-\psi\right)\zetadd\right]\ccu\vv dxdt \\
&\;\;\;\, +\frac{1}{\delta} \int^h_0\int_{\left[0<d_\delta <\delta\right]}\left[\left(\psi_C-\psi\right)\ccu\right](t,s,\nu,\omega)\nabla d_\eps\cdot\nabla F^{-1}\vv(t,s,\nu,\omega)J_Fdx dt
\end{align*}

The integral for the initial time trace follows again because $\zetadd\rightarrow 0$ a.e.
\hfill \bbox
\begin{remark}
	If $\vv$ is assumed smooth, then both the time derivative and advection terms could be handled by coarea arguments that circumvent the need for $\psi_C-\psi = O(\delta)$. However, the diffusion term could not since $\nabla\ccu$ is not sufficiently regular to have traces on hypersurfaces amenable to the divergence theorem.
\end{remark}
\begin{corollary}\label{cr:testvphi}
For a given smooth $\psi\in\mathcal{C}^\infty\left(\bar{\Omega}_h\right)$, the choices 
\[
\psi_C = \psi(t,s,0,0),\psi(t,s,\nu,0),\psi(t,s,0,\omega),\psi(t,s,\nu,\omega)
\]
all satisfy $\left|\psi - \psi_C\right| = O(\delta)$ over $\mathcal{N}_{\delta,\delta}(t)$. Accordingly, for these choices (\ref{eq:wklim}) holds with $\varphi$ replaced by $\psi$ in the terms $\left(\ldots\right)_{\mathrm{Vol}}$ and replaced by $\psi_C$ in the terms $\left(\ldots\right)_{\mathrm{Line}}$. 
\end{corollary}
\emph{Proof:} Since the image of $\mathcal{F}$ is compact and contained in $\Omega$, $\nabla\psi$ takes a global maximum which may be used in conjunction with a line integral over the convex domain of $\mathcal{F}$ to bound the difference. We thus satisfy the conditions of \textbf{Lemma \ref{lm:testvphi}}. \hfill \bbox
\subsection{Concluding}
\begin{proposition}
	The weak limits $\xi_\nu,\xi_\omega\in L^2\left(\left[0,1\right]_h\right)$ are given by (\ref{eq:xiwklim}).
\end{proposition}
\emph{Proof:} Take the weak formulation supplied by \textbf{Corollary \ref{cr:testvphi}} and subtract the result for the pair $\big(\psi,\psi(t,s,0,0)\big)$ from the result for the pair $\big(\psi,\psi(t,s,\nu,0)\big)$. Since both pairs produce the same terms in the volume, we need only examine the terms at the curve. After canceling we are left with
\begin{align*}
0 &= \int_{\left[0,1\right]_h}\bigg(\partial_\nu\psi\vec{e}_2\cdot\nabla F^{-1}\partial_tF\ccuC - \partial_\nu\psi\vec{e}_2\cdot\nabla F^{-1}\vv_C\ccuC+\partial_\nu\psi k_n\xi_\nu\bigg)\JG dsdt \\
&= \int_{\left[0,1\right]_h}\partial_\nu\psi\left(\vec{e}_2\cdot\nabla F^{-1}\left(\partial_t F - \vv_C\right)\ccuC + k_n\xi_v\right)\JG dsdt
\end{align*}
From the arbitrariness of $\psi$ in the volume, $\partial_\nu\psi$ can recover all smooth testing functions on $\left[0,1\right]_h$ in the variables $(t,s)$ which vanish at time $t = h$. We conclude the integrand is identically 0, and thus $\xi_\nu$ must be
\[
\xi_\nu(t,s) = \frac{\left[\vec{e}_2\cdot\nabla F^{-1}\left(\vv_C-\partial_tF \right)\right](t,s,0,0)\ccuC}{k_n(s,0,0)}
\]
The function $\xi_\omega$ may be derived analogously.
\[
\xi_\omega(t,s) =  \frac{\left[\vec{e}_3\cdot\nabla F^{-1}\left(\vv_C-\partial_tF \right)\right](t,s,0,0)\ccuC}{k_n(s,0,0)}
\]
\hfill \bbox
\section{Appendix}
\subsection{Coordinate Transformations and Geometric Properties}
In this appendix we collect several elementary but helpful observations regarding changes of variable and the geometry used in this investigation. The quantities in \textbf{Lemma \ref{lm:DF}} and \textbf{Lemma \ref{lm:DFinv}} are required for administration of chain-rule.
\begin{lemma}Jacobian Matrices\label{lm:DF}\\
It holds that
\begin{align*}
D\mathcal{F} &= 
\begin{pmatrix}
1 & \vec{0}\\
\partial_t F & \nabla F\left(=D_{s,\nu,\omega}F\right)
\end{pmatrix}\\
\nabla F &= \begin{pmatrix} \partial_s\Gamma + \nu\partial_s\vn + \omega\partial_s\vb & \vn & \vb\end {pmatrix}\\
\nabla F^\mathrm{t}\nabla F &=
\begin{pmatrix}
\left\|\partial_s\Gamma + \nu\partial_s\vn+\omega\partial_s\vb \right\|^2 & \omega\left\langle \partial_s\vb,\vn \right\rangle & \nu\left\langle \partial_s\vn,\vb\right\rangle\\
\omega\left\langle \partial_s\vb,\vn \right\rangle & 1 & 0\\
\nu\left\langle \partial_s\vn,\vb\right\rangle & 0 & 1
\end{pmatrix}\\
J_{\mathcal{F}} &:= \mathrm{det}\, D\mathcal{F} = \mathrm{det}\, \nabla F\\
&=\sqrt{\left|\partial_s\Gamma + \nu\partial_s\vn+\omega\partial_s\vb \right|^2 - \nu^2\left\langle \partial_s\vn,\vb\right\rangle - \omega^2\left\langle \partial_s\vb,\vn\right\rangle^2}=: J_{F}
\end{align*}
\end{lemma}
\emph{Proof:} Use (\ref{eq:DF}) and elementary manipulations.\hfill \bbox
\begin{lemma}Inverses of Jacobian Matrices\label{lm:DFinv}\\
It holds that
\begin{align*}
D\mathcal{F}&^{-1} = 
\begin{pmatrix}
1 & \vec{0}\\
-\nabla F^{-1}\partial_t F & \nabla F^{-1}
\end{pmatrix}\\
\left(\nabla F^{\mathrm{t}}\nabla F\right)&^{-1}=\\ & \frac{\begin{pmatrix}
	1 & -\omega\left\langle \partial_s\vb,\vn\right\rangle & -\nu\left\langle \partial_s \vn,\vb \right\rangle \\
	-\omega\left\langle \partial_s\vb,\vn\right\rangle & \left|\partial_s\Gamma + \nu\partial_s\vn+\omega\partial_s\vb \right|^2 - \nu^2\left\langle \partial_s\vn,\vb\right\rangle & \nu\omega\left\langle\partial_s\vn,\vb\right\rangle\left\langle\partial_s\vb,\vn\right\rangle\\
	-\nu\left\langle \partial_s \vn,\vb \right\rangle & \nu\omega\left\langle\partial_s\vn,\vb\right\rangle\left\langle\partial_s\vb,\vn\right\rangle & \left|\partial_s\Gamma + \nu\partial_s\vn+\omega\partial_s\vb \right|^2 - \omega^2\left\langle \partial_s\vb,\vn\right\rangle^2
	\end{pmatrix}}{\left|\partial_s\Gamma + \nu\partial_s\vn+\omega\partial_s\vb \right|^2 - \nu^2\left\langle \partial_s\vn,\vb\right\rangle - \omega^2\left\langle \partial_s\vb,\vn\right\rangle^2}\\
\nabla F^{-1} &= \left(\nabla F^{\mathrm{t}}\nabla F\right)^{-1} \nabla F^{\mathrm{t}}\\
&= \left(\nabla F^{\mathrm{t}}\nabla F\right)^{-1} 
\begin{pmatrix}
\partial_s\Gamma^{\mathrm{t}} + \nu\partial_s\vn^{\mathrm{t}} + \omega\partial_s\vb^{\mathrm{t}}\\
\vn^{\mathrm{t}}\\
\vb^{\mathrm{t}}
\end{pmatrix}
\end{align*}
\end{lemma}
\emph{Proof:} Use \textbf{Lemma \ref{lm:DF}} and elementary manipulations. \hfill \bbox 

The distance function $d_\eps(x)$ may be computed explicitly due to the simplicity of $\mathcal{C}_\eps$.

\begin{lemma}The distance function \label{lm:deps}\\
	For $x=\left(s,\nu,\omega\right)$
	\begin{equation*}
	d_\eps(x) = \sqrt{s^2_-+(s-1)^2_++\left(\sqrt{\nu^2+\omega^2}-\eps\right)^2_+}
	\end{equation*}
\end{lemma}
\emph{Proof:} We may decompose the distance function and write
\[
\mathrm{min}_{y\in\mathcal{C}_\eps} \left\|x-y\right\|^2 = \text{dist}\big(s,\left[0,1\right]\big)^2 + \text{dist}\big(\left(\nu,\omega\right),D_\eps\big)^2
\]
\hfill \bbox
\begin{lemma}Measuring $\left[\mathcal{N}_{\eps,\delta}\setminus\mathcal{N}_{\eps,0}\right](t)$ \label{lm:gap}\\
There exists a constant $C = C(\eps_0)$ s.t.
\begin{align*}
\mu\bigg(\left[0<d_\eps<\delta\right]\bigg) &\leq C(\eps_0)\delta\\
\mu\bigg(\left[\Ned\setminus\mathcal{N}_{\eps,0}\right](t)\bigg) &\leq\left( \esup_{t\leq T}\left\|J_F\right\|_{\infty}(t)\right) C(\eps_0)\delta
\end{align*} 	
\end{lemma}
\emph{Proof:} Begin by applying the coarea formula.
\[
\mu\bigg(\left[0<d_\eps<\delta\right]\bigg) = \int_{\left[0<d_\eps<\delta\right]}1 dx = \int^\delta_0\int_{\left[d_\eps = \eta \right]}\frac{1}{\left|\nabla d_\eps\right|}d\sigma d\eta
\]
From \textbf{Lemma \ref{lm:deps}}, one verifies that $\left|\nabla d_\eps(x)\right| = 1$ when $d_\eps(x)\neq 0$. From the geometry of $\mathcal{C}_\eps$, we have $\exists C(\eps_0)$ s.t. $\forall \tau\in\left(0,\delta\right)$
\[
\int_{\left[d_\eps = \tau\right]}d\sigma \leq C(\eps_0)
\]
This concludes the first assertion.  The second follows from a simple change of variables.
\[
\int_{\left[\Ned\setminus\mathcal{N}_{\eps,0}\right](t)}dx = \int_{\left[0<d_\eps < \delta\right]}J_F dx\leq \left( \esup_{t\leq T}\left\|J_F\right\|_{\infty}(t)\right) \mu\bigg(\left[0<d_\eps<\delta\right]\bigg)
\]
\hfill \bbox
\subsection{The Distributional Limit of $\aed(t;x)$}
In this appendix, we prove a simple lemma regarding the weak convergence of $\aed(t;x)$ as distributions on the space of continuous functions.
\begin{lemma} The weak limit of $\aed$\label{lm:aedwklim}\\
Let $f\in \mathcal{C}(\bar{\Omega})$ and $\eps_i\rightarrow 0$. Then $\exists \kappa_i\left(>0\right)\rightarrow 0$, so that if $\delta_i<\kappa_i$ it holds
\[
\lim_i\int_\Omega a_{\eps_i,\delta_i}(t;\cdot)fdx = \int_\Omega f dx + \pi\eps^2_0\int^1_0\left[\left(f\circ\Gamma\right)\JG\right](t;s) ds
\]
\end{lemma}
\textit{Proof:} Begin by writing
\[
\int_\Omega a_{\eps_i,\delta_i} f dx = \int_\Omega \left(1-\zeta_{\eps_i,\delta_i}\right)f dx + \int_\Omega \frac{\eps^2_0}{\eps^2_i}\zeta_{\eps_i,\delta_i} f dx
\]
As $\eps_i\rightarrow 0$, $\chi_{\mathcal{N}_{\eps_i,0}(t)}\rightarrow 0$ a.e. and so does strongly as well. Since $\left|\zeta_{\eps_i,\delta_i}-\chi_{\mathcal{N}_{\eps_i,0}(t)}\right| \leq 2$ and \textbf{Lemma \ref{lm:gap}} shows the set where they differ has measure $O(\delta_i)$, we have 
\[
\left|\int_\Omega\left(1-\zeta_{\eps_i,\delta_i}\right)f dx - \int_\Omega f dx \right| \leq 2C(\mathrm{data})\left\|f\right\|_\infty\delta_i + \left|\int_\Omega\left(1-\chi_{\mathcal{N}_{\eps_i,0}(t)}\right)fdx - \int_\Omega f dx\right|
\]
This proves the result for the first term. We proceed as follows to handle the second.
\begin{align*}
\int_\Omega\frac{\eps^2_0}{\eps^2_i}\zeta_{\eps_i,\delta_i} fdx &= \int_{\mathcal{N}_{\eps_0}(t)}\frac{\eps^2_0}{\eps^2_i}\zeta_{\eps_i,\delta_i} \left(f\circ F\right) dx = \int_{\left[-\eps_0,1+\eps_0\right]\times D_{\eps_0}}\frac{\eps^2_0}{\eps^2_i}\left(\chi_{\delta_i}\circ d_{\eps_i}\right)\left(f\circ F\right)J_{F}dx\\
&= \int_{\mathcal{C}_\eps}\frac{\eps^2_0}{\eps^2_i}\left(f\circ F\right) J_{F}dx + \int_{\left[0<d_{\eps_i} <{\delta_i}\right]}\frac{\eps^2_0}{\eps^2_i}\left(\chi_{\delta_i}\circ d_{\eps_i}\right)\left(f\circ F\right)J_{F}dx\\
&= \pi\eps^2_0\int^1_0\dashint_{D_{\eps}}\left(f\circ F\right) J_F dx + \int_{\left[0<d_\eps <\delta\right]}\frac{\eps^2_0}{\eps^2_i}\left(\chi_{\delta_i}\circ d_{\eps_i}\right)\left(f\circ F\right)J_{F}dx
\end{align*}
From \textbf{Lemma \ref{lm:gap}}, we may write
\[
\int_\Omega\frac{\eps^2_0}{\eps^2_i}\zeta_{\eps_i,\delta_i} fdx = \pi\eps^2_0\int^1_0\dashint_{D_{\eps_i}}\left(f\circ F\right) J_F dx + C(\mathrm{data})\left\|J_F\right\|_\infty \left\|f\right\|_\infty\frac{\delta_i}{\eps^2_i}
\]
The result is concluded using $f$'s and $J_F$'s continuity. Note that \textbf{Lemma \ref{lm:DF}} implies that $J_F\rightarrow\JG$. \hfill \bbox
\begin{remark}\label{rm:Diwklim}
	We observe for concreteness that we may take $\kappa_i = \eps^3_i$.
\end{remark}

\end{document}